\documentclass{article}

\usepackage{amssymb}
\usepackage{amsmath}

\let\<\langle
\let\>\rangle
\font\bigcmsy=cmsy10.pk scaled 2000
\def\bigtimes{\mathop{\,\vrule width0pt depth2pt height8pt
            \smash{\lower2pt\hbox{\bigcmsy\char'002}}\,}\limits}

\begin{document}

\begin{center}
\Large{\textbf{Order separability of free product of free groups with cyclic amalgamation.}}
\end{center}
\begin{center}

\textbf{Vladimir V. Yedynak}

\end{center}

\begin{abstract}

This paper is devoted to the proof of the property of order separability for free product of free groups with maximal cyclic amalgamated subgroups.


\end{abstract}

\section{Introduction.}

Definition 1. \textsl{Consider the group $G$ and the natural number $n$, which is greater than 1. Then we say that the group $G$ is $n$-order separable if for each elements $g_1,..., g_n$ of $G$ such that $g_i$ is conjugate to neither $g_j$ not $g_j^{-1}, 1\leqslant i< j\leqslant n,$ there exists a homomorphism $\varphi$ of the group $G$ onto a finite group such that the orders of the elements $\varphi(g_1),..., \varphi(g_n)$ are pairwise different}.

It was proved in [1] that free groups are 2-order separable. Later in [2] this result was generalized by Daniel T. Wise who showed that free groups are omnipotent. Different properties of groups connected with order separabilities such as potency, quasipotency, omnipotency, 2-order separability enable to investigate residual properties of groups. For example in [3] Graham A. Niblo investigated the existence of \textsl{regular quotients} at $\{U, V\}$ in order to obtain classes of HNN extensions of free groups with cyclic connected subgroups which are subgroup separable. Here $U$ and $V$ are the elements of a free group $F$ and $F$ is said to have regular quotients at $\{U, V\}$ if there exists a natural number $m$ such that for each natural $n$ there exists a homomorphism of $F$ onto a finite group such that the images of $U$ and $V$ both have order $nm$. There in after we prove the generalization of this result spreading it to a finite set of elements of a free group for to investigate the order separability of free products with amalgamation.

The property of order separability was also investigate in [4] where it was proved that 2-order separability is inherited by free products. The aim of this work is to prove the following theorem.

\textbf{Theoren 1.} \textsl{Let $F_1$ and $F_2$ be a free groups. Consider elements $a\in F_1\setminus\{1\}, b\in
F_2\setminus\{1\}$ such that $a$ and $b$ are not proper powers that is $a\neq w^k, |k|>1, w\in F_1, b\neq z^l, |l|>1, z\in F_2$. Then the group $F_1{\ast\atop a=b}F_2$ is 2-order separable.}

\section{Auxiliary facts from the graph theory.}

In this section we give basic facts of the graph theory that will be useful in what follows.

Definition 1. \textsl{Graph} $\Gamma$ is a pair of sets $E,
V$, for which three maps derived: $\alpha: E\rightarrow V,
\omega: E\rightarrow V, \eta: E\rightarrow E$ such that the following conditions are true: $\eta^2=$ id, $\alpha\circ\eta=\omega,
\omega\circ\eta=\alpha$, the map $\eta$ is the bijection and does not possess fixed points.

The sets $E$ and $V$ are called the edges and vertices of the graph $\Gamma$.
The map $\alpha$ corresponds for each edge of $\Gamma$ its begin point, while the map
$\omega$ corresponds for each edge its end point.
The map $\eta$ corresponds for each edge $e$ the edge $f$ which is inverse to $e$. Edges $e$ and $f$ are called mutually inverse.

Definition 2. The graph is called \textsl{oriented}, if from each pair of its mutually inverse edges one of them is fixed, it is called positively oriented. Edges which are inverse to positively oriented are called negatively oriented.

Definition 3. The set of edges $\{e_i| i\in I\}$, where either
$I=\{1,2,..., n\}$ or $I$ coincide with the set of natural numbers,
of the graph $\Gamma$ is called a \textsl{path} $S=e_1...e_n...$, if
$\alpha(e_i)=\omega(e_{i-1})$ for each $i$ such that $i-1, i\in I$.
The vertex $\alpha(e_1)$ is called the beginning of $S$ and will be denoted by $\alpha(S)$. If $I=\{1, 2,..., n\}$ then the vertex $\omega(e_n)$ is called the end of the path $S:$ $ \omega(S)$. In this case we shall also say that the path $S$ is finite.

Definition 4. The path $S=\{e_i| i\in I\}$ is called the \textsl{path without returns}, if edges $e_i, e_{i+1}$ both are not mutually inverse for each $i, i+1\in I$. In what follows we shall consider only such paths.

Definition 5. The finite path is called a \textsl{cycle} if its beginning and end coincide.

Definition 6. Consider two vertices $p$ and $q$ of the graph $\Gamma$ which belong to one component of $\Gamma$. Then the
\textsl{distance} between $p$ and $q$ is the minimal length of the finite path which connects $p$ and $q$. The length of the finite path is the number of edges which compose this path.

Definition 7. Consider the graph $\Gamma$, the cycle $S=e_1...e_n$ in it and the nonnegative integer number $l$. Wa say that $S$ does not have \textsl{$l$-near vertices} if for each $i, j, 1\leqslant i<j\leqslant n$ the distance between the vertices
$\alpha(e_i), \alpha(e_j)$ is greater or equal than min$(|i-j|, n-|i-j|, l+1)$.

Note that when $l=0$ this condition means that $S$ does not have self-intersections. In this case we shall say that $S$ is simple. When $l=1$ we shall say that $S$ does not have near vertices. Remark also that if $S$ does not have $l$-near vertices then $S$ does not have $k$-near vertices for each nonnegative integer $k$ which is less than $l$.

\section{Action graphs of groups.}

Consider a group $G$ acting on the set $X$. Fix a generating set of
$G: \{g_i| i\in I\}$. Consider the oriented graph $\Gamma$ in the following way. Elements of
$X$ are identified with vertices of $\Gamma$. For each vertex $p$ and for each $i\in I$ there exist vertices $q$ and $r$ such that $p\circ g_i=q, p\circ g_i^{-1}=r$. Connect vertices $p$ and $q$ by the pair of mutually inverse edges. Furthermore the edge which goes away from $p$ into $q$ is regarded positively oriented and has label $g_i$. Connect also the vertices $p$ and $r$ by the pair of mutually inverse edges and the edge which goes away from $r$ into $p$ is positively oriented and has label $g_i$.

Definition 1. The derived graph is called the \textsl{action graph} of the group $G$ with respect to the generating set $\{g_i|i\in I\}$. The label of the positively oriented edge $f$ will be denoted as Lab$(f)$. We shall talk about action graphs of some group $G$ without references on its set of generators if it is fixes.

Definition 2. Consider the path $S=e_1...e_n$ in the action graph $\Gamma$
of the group $G=\<g_i|i\in I\>$. Then the \textsl{label of the path $S$}
is the element of the group $G$ which is equal to $\prod_{i=1}^n$Lab$(e_i)'$,
where Lab$(e_i)'=$ Lab$(e_i)$, if the edge $e_i$ is positively oriented. And if the edge $e_i$ is negatively oriented then
Lab$(e_i)'=$ Lab$(\eta(e_i))^{-1}$.

Definition 3. Consider the element $u$ of the group $G$ which does not equal to unit. Let $\Gamma$ be the action graph of the group $G=\<g_i|i\in I\>$. Fix the vertex $p$ in $\Gamma$. Then the \textsl{$u$-cycle $T$} of the graph $\Gamma$ going from $p$ is the set of paths $\{S_i=\{e_j^i| j\in J_i\}|i\in J\}$ satisfying the following conditions:

1) $\alpha(S_l)=p$;

2) there exists a one to one correspondence between the paths $S_l$ and the presentations of $u$ in generators $\{g_i|i\in I\}: u=g_{i_1}^{\varepsilon_1}...g_{i_k}^{\varepsilon_k}, $
min$(|i_j-i_{j+1}|, \varepsilon_j+\varepsilon_{j+1})>0,
\varepsilon_i\in\{-1; 1\}, 1\leqslant j\leqslant k, 1\leqslant
i\leqslant k$. Besides if $\varepsilon_j=1$ then edge $e_{kn+j}^l$ is positively oriented and its label equals $g_{i_j}$ for each natural $n$; if $\varepsilon_j=-1$ then $e_{kn+j}^l$ is negatively oriented and the label of the edge $\eta(e_{i_j})$ equals $g_{i_j}$ for each natural $n$; in case the path $S_l$ is finite and is composed of $r$ edges we consider that indices are modulo $r$;

3) there is no closed subpath $K$ in the path $S_l$ which differs from $S_l$ and satisfies the conditions 1), 2);

4) the path $S_l$ is the maximal path on the entry which satisfies the conditions 1)-3);

Paths of the set $\{S_i|i\in J\}$ are called the representatives of the $u$-cycle $T$.

It is easy to notice that the representatives of $T$ contains elements of the orbit of the vertex $p$ at the action of positive powers of $u$ on $p$. And if the graph is finite then there exists a one to one correspondence between $u$-cycles of $\Gamma$ and orbits at the action of the subgroup $\<u\>$ on the set of vertices of $\Gamma$.

Note that if some representative of the $u$-cycle $T$ contains only a finite set of edges then it is closed and all the rest representatives of $T$ are also closed and the label of each such representative is equal to $u^k$ for some $k$.

Definition 4. Suppose that some representative of the $u$-cycle $T$ is composed from the finite set of edges and has the label  $u^k$. Then we shall say that the \textsl{length} of the $u$-cycle $T$ is equal to $k$.

Let $G$ be the group generated by the set $\{g_i| i\in I\}$. Consider the epimorphism $\varphi$ from the free group $F$ with the basis $\{x_i| i\in I\}$ onto the group $G$ such that the image of the element $x_i$ is equal to the element $g_i$ for each $i\in I$. Consider that the kernel of $\varphi$ is the normal closure of elements $\{R_j|
j\in J\}$. Then the oriented graph $\Gamma$ whose positively oriented edges are labelled by the elements from $\{g_i| i\in
I\}$ is the action graph of the group $G$ with the generating set $\{g_i| i\in I\}$ if and only if the following conditions are true (we suppose that in $\Gamma$ each label $x_i$ is identified with $g_i$ for each $i\in I$):

1) for each vertex $p$ of $\Gamma$ and for each $i\in I$ there exist exactly one positively oriented edges $e$ and $f$, labelled by $x_i$ such that $f$ goes into $p$ and $e$ goes away from $p$;

2) for each vertex $p$ from $\Gamma$ and for each $j\in J$ the $R_j$-cycle going from $p$ has length one;

It is evident also that there is a correspondence between the actions of the group $G$ and the oriented graphs satisfying the properties 1), 2).

Definition 5. Consider the group $G$ with the set of generators $\{g_i|i\in I\}$ and the action graph $\Gamma$ of the group $G=\<g_i|i\in I\>$. Then the homomorphism of $G$ onto the permutation group on the set of vertices of $\Gamma$ will be denoted as \textsl{$\varphi_{\Gamma}$}.

Consider now the action graph $\Gamma$ of the group $G=\<g_i|i\in I\>$ which possesses the finite number of vertices. It is true that for each $u$-cycle $T$ from $\Gamma$ each representative of $T$ consists of the finite set of edges. Hence the order of the element $\varphi_{\Gamma}(u)$ is finite for each $u\in G$ and equals the least common multiple of lengths of all $u$-cycles in $\Gamma$. And if in the graph $\Gamma$ there exists the $u$-cycle whose length is equal to $t$ then $\mid\varphi_{\Gamma}(u)\mid\geqslant t$.

Working with action graphs we shall append and delete oriented edges with labels. We shall always suppose that after deleting the edge with label the edge which is inverse to it is also deleted. And if we append an edge with label the corresponding inverse edge would have been appended.

\section{Auxiliary notations and lemmas.}

Definition. \textsl{A group $G$ is called subgroup separable if for each element $g$ of $G$ and for each finitely generated subgroup $K$ of $G$ such that $g\notin K$ there exists a finite quotient $G/N$ such that $gN\notin KN$. In other words $G$ is subgroup separable if each finitely generated subgroup of $G$ is a closed set in the profinite topology  }

Consider that $\Gamma$ is the action graph of the group $A{\ast\atop
C}B$ with the generating set $\{A\cup B\}$ such that the subgroups $A$ and $B$ act freely on vertices of $\Gamma$. Then it is equivalent to the realization of the following conditions for the graph $\Gamma$:

(1) $\Gamma$ is the oriented graph whose positively oriented edges are labelled by the elements of $A\cup B$ such that for each vertex $p$ of $\Gamma$ and for each $d\in A\cup B$ there exist exactly one edge with label $d$ going from $p$ and exactly one edge with label $d$ going into $p$

(2) for each vertex $p$ from $\Gamma$ we define the subgraph $A(p)$ as the maximal connected subgraph containing $p$ whose positively oriented edges are labelled by the elements from $A$; the graph $A(p)$ is the Cayley graph of the group $A$ with the generating set $\{A\}$; the subgraph $B(p)$ is defined analogically;

Derive the following notation. If $G$ is an arbitrary group and $V$ its generating set then $Cay(G, V)$ is the Cayley graph of the group $G$ with respect to its generating set $V$.

In what follows we suppose that all concerned action graphs of free product with amalgamation are so that free factors act freely on vertices of action graphs. For the free product $G$ the generating set of $G$ will be always the union of its free factors.

The following theorem was proved in [5].

\textbf{Theorem 2.} \textsl{Class of subgroup separable groups is closed with respect to the operation of free product.}

From the theorem 2 we may deduce that free groups are subgroup separable.

In [2] the following well known result was proved geometrically.

\textbf{Theorem 3.} If a group $G$ contains a finite index subgroup which is subgroup separable then $G$ is also subgroup separable.

Definition 1. A group $G$ is called potent if for each element $g$ of $G$ and for each natural number $n$ there exists a normal divisor $H_n$ of $G$ of finite index such that $\<g\>\cap H_n=\<g^n\>$.

The following theorem was proved in [6].

\textbf{Theorem 4.} \textsl{Any absolutely free group is potent.}

Therein after we shall identify $u$-cycles and the corresponding representatives of $u$-cycles for action graphs of free groups with respect to the set of generators coinciding with the fixed basis.

In what follows we shall use the following notation. Consider the action graph $\Gamma$ and a subset $R$ of $\Gamma$. Having a set of copies of $\Gamma: \Delta_1,..., \Delta_n$ we shall suppose that $R^i$ is the subset of $\Delta_i$ corresponding to $R$.

\textbf{Lemma 1.} \textsl{Let $F$ be a free group with a fixed basis $Z$, $u_1,..., u_k, v\in F$ and for each $i$ elements $v$ and $u_i$ does not belong to conjugate cyclic subgroups. Then for each quit great prime number $p$ and for each natural number $N$ there exists a homomorphism $\varphi$ of $F$ onto a finite $p$-group such that $\mid\varphi(u_1)\mid=...=\mid\varphi(u_n)\mid>\mid\varphi(v)\mid>1$,
and $\mid\varphi(u_1)\mid>N$.}

\textbf{Proof.}

Consider the prime number $p$ such that for each $i$ \ $u_i\neq u_i'^p, v\neq v'^p$. Hence if $u_i, u_j$ belong to conjugate cyclic subgroups then their images have equal orders after any homomorphism of $F$ onto a finite $p$-group. So we shall consider that the elements $u_1,..., u_k$ are pairwise incommensurable. Put $u_i=u_i'^{m_i}, v=v'^l$, where $u_i', v'$ are not proper powers. According to the conditions on the prime $p$ we deduce that if there exists the required homomorphism for elements $u_1',..., u_n', v'$ then there also exists the required homomorphism for elements $u_1,..., u_n, v$. So we may consider that the elements $u_1,..., u_n, v$ are not proper powers, particularly subgroups $\<u_1\>,..., \<u_k\>, \<v\>$ are $p'$-isolated.

Consider the following case. Consider that for some $m, 1\leqslant m\leqslant n-1$, there exists a homomorphism $\varphi_m$ of $F$ onto a finite $p$-group such that in the graph $Cay(\varphi_m(G); \varphi_m(Z))$ all $v$- and $u_j$-cycles are simple, $1\leqslant j\leqslant n$. In addition for each $i, m+1\leqslant i\leqslant k, \mid\varphi_m(u_1)\mid=...=\mid\varphi_m(u_m)\mid<\mid\varphi_m(u_i)\mid$. These conditions are true for some $m\geqslant1$. According to [7] there exists a homomorphism $\varphi_1$ of $F$ onto a finite $p$-group such that $\varphi_1(u_1^j),..., \varphi_1(u_k^j), \varphi_1(v^j)\neq1, 1\leqslant j\leqslant N$. Particularly $\mid\varphi_1(u_i)\mid>N, 1\leqslant i\leqslant k$. Since free groups are residually $p$-finite with respect to the entry into the $p'$-isolated cyclic subgroups and are residually $p$-finite we may consider that in the Cayley graph of the group $\varphi_1(F)$ with generators $\varphi_1(Z)$ all
$u_i$- and $v$-cycled are simple, $1\leqslant i\leqslant
n$. After this the elements whose images after $\varphi_1$ have minimal orders may taken in the capacity of $u_1,..., u_m$.

Consider the Cayley graph $\Delta=Cay(\varphi_m(G); \{\varphi_m(Z)\})$. Consider the inductive process of constructing of action graphs of the group $F=\<Z\>$ $\Gamma_{k}, 1\leqslant k\leqslant w$ (the number $w$ will be chosen later) satisfying the following properties.

1) the length of each $v$- and $u_i$-cycle divides the length of the maximal $v$-,
correspondingly $u_i$-cycle, $1\leqslant i\leqslant n$.

2) all $v$- and $u_i$-cycles are simple, $1\leqslant i\leqslant n$.

3) for $k, 1\leqslant k\leqslant w$ in the graph $\Gamma_{k}$ there exists the path $S_k$ of length $k$ which is contained in the maximal $u_1$-cycle and in all maximal $u_l$-cycles, $m+1\leqslant l\leqslant n$.

Construct the graph $\Gamma_{1}$ first. Consider $p$ copies of the graph $\Delta: \Delta_1,..., \Delta_p$. In the graph $\Delta$ fix some $u_1$-cycle $S=e_1...e_s$. Delete edges $e_1^1,..., e_1^p$ from the graphs $\Delta_1,..., \Delta_p$. For each $i, 1\leqslant i\leqslant p$, if the edge $e_1$ is positively oriented then we append a new positively oriented edge $f_i$ going from the vertex $\alpha(e_1^i)$ into the vertex $\omega(e_1^{i+1})$. If the edge $e_1$ is negatively oriented then the derived edge $f_i$ goes from the vertex $\omega(e_1^{i+1})$ and goes into the vertex $\alpha(e_1^i)$. In both cases the label of the edge $f_i$ coincide with the label of the positively oriented edge from the pair of edges $e_1, \eta(e_1)$ (all upper indices are modulo $p$). We shall denote the obtained action graph of the group $F$ as $\Gamma_{1}'$. Conditions 1), 2) are true for $\Gamma_{1}'$ by the construction. It is obvious also that there exists a maximal $u_1$-cycle passing through the edge $f_1$. If there exists $l, m+1\leqslant l\leqslant n$ such that there is no $u_l$-cycle containing the edge $f_1$ then we consider that the process has been finished and $w=0$ and we also fix the number $l$ such that there exists the maximal $u_l$-cycle which does not pass through the edge $f_1, l>m$. Otherwise we put $\Gamma_1=\Gamma_1', S_1=f_1$.

Suppose we have the graph $\Gamma_{k}, k\geqslant1$ satisfying properties 1), 2), 3). Fix the maximal $u_1$-cycle $T=q_1...q_d$ passing through $S_k$. By the inductive supposition we may consider that $S_k=q_1...q_k$. Put $q=\omega(q_k)$. Consider $p$ copies of the graph $\Gamma_{k}: \Delta_{k}^1,..., \Delta_{k}^p$. Delete edges $q_{k+1}^1,..., q_{k+1}^1$ from the graphs $\Delta_{k}^1,..., \Delta_{k}^p$. For each $i, 1\leqslant i\leqslant p$ derive the new positively oriented edge $g_i$. If the edge $q_{k+1}$ is positively oriented then the edge $g_i$ goes away from $\alpha(q_{k+1}^i)$ and goes into $\omega(q_{k+1}^{i+1})$. If $q_{k+1}$ is negatively oriented then $g_i$ goes away from $\omega(q_{k+1}^{i+1})$ and goes into $\alpha(q_{k+1}^i)$. The label of $g_i$ coincides with then label of the positively oriented edge from the pair $q_{k+1}, \eta(q_{k+1})$ (all upper indices are modulo $p$). For the obtained action graph $\Gamma_{k+1}'$ properties 1), 2) are held by the construction. Furthermore the path $q_1^1...q_k^1\eta'(g_1)$ is contained in the maximal $u_1$-cycle where $\eta'(g_1)$ is the edge from $\{g_1, \eta(g_1)\}$ whose beginning coincides with the end of the edge $q_k^1$. If this edge is contained in all maximal $u_l$-cycles, $m+1\leqslant l\leqslant n$, then according to property 2) for $\Gamma_k$ these cycles pass through $q^1_1...q^1_k\eta'(g_1)$ also and in this case we denote this path as $S_{k+1}$ and put $\Gamma_{k+1}=\Gamma_{k+1}'$. Otherwise we consider that the process has been finished and $w=k$ and we fix the number $h, m+1\leqslant h\leqslant n$ such that there exists the maximal $u_h$-cycle which does not pass through the path $q_1^1...q_k^1\eta'(g_1)$. Since the elements $u_1,..., u_n$ are pairwise incommensurable the process will be finished with some finite number $w$ and we shall obtain eventually the element $u_h, m+1\leqslant h\leqslant n,$ such that in the graph $\Gamma_{w+1}'$ not all maximal $u_h$-cycles go through the path of length $w$.

Now we shall work with the element $u_h$. Fix $i, 1\leqslant i\leqslant m$. Consider the inductive process of building the action graphs $\Lambda_{i1},..., \Lambda_{i, r_i}$ such that the graphs satisfy the following properties.

1) the length of each $v$- and $u_s$-cycle divide the length of the maximal $v$-,
correspondingly $u_s$-cycle, $1\leqslant s\leqslant n$.

2) all $v$- and $u_s$-cycles are simple, $1\leqslant s\leqslant n$.

3) for $k, 1\leqslant k\leqslant r_i-1$ in the graph $\Lambda_{ik}$ there exists the path $S_i^k$ of length $k$ which is contained in some maximal $u_i$-ccyle and in all maximal $u_h$-cycles.

The aim of the first inductive process was just to find the element $u_h$ which is one of the first for which the variant of the condition 3) is not held for constructible graphs. The second process does not differ from the previous inductive process except for the modified condition 3). Besides on the last step we consider not $p$ copies of the graph $\Lambda_{i, r_i-1}$ but $p^{q_i}$ copies. The numbers $q_i$ are so that for each $i, j, 1\leqslant i, k\leqslant m$ the following equality is true: $r_i+q_i=r_j+q_j$. Let's estimate the lengths of maximal $u_i$-cycles, $1\leqslant i\leqslant n,$ in the graph $\Lambda=\bigsqcup_{i=1}^m\Lambda_{i, r_i}$. By the construction and due to equality $|\varphi_m(u_1)|=...=|\varphi_m(u_m)|$ we deduce that in the graph $\Lambda_{i, r_i}$ the length of each $u_s$-cycle, $1\leqslant s\leqslant m, s\neq i,$ is the power of $p$ and is not greater than the length of the maximal $u_i$-cycle in the graph $\Lambda_{i, r_i}$. Hence using the condition on numbers $q_l$ we conclude that $|\varphi_{\Lambda}(u_1)|=...=|\varphi_{\Lambda}(u_m)|$. Let's show that $|\varphi_m(u_h)|/|\varphi_m(u_1)|<|\varphi_{\Lambda}(u_h)|/|\varphi_{\Lambda}(u_1)|$. Every time we construct the graph $\Lambda_{i, k+1}$ from the graph $\Lambda_{i, k}$ when $k\neq r_i-1$ the lengths of $u_i$- and of $u_h$-cycles increase in $p$ times. When we construct the graph $\Lambda_{i, r_i}$ from $\Lambda_{i, r_i-1}$ the length of the maximal $u_i$-cycle increases in $p^{q_i}$ times and the length of the maximal $u_h$-cycle increases no more than in $p^{q_i-1}$ times because for each vertex $q$ of the graph $\Lambda_{i, r_i-1}$ which is the beginning of a maximal in the graph $\Lambda_{i,
r_i-1}$ $u_h$-cycle the corresponding $u_h$-cycle of the graph $\Lambda_{i, r_i}$ going from the vertex $q_i, 1\leqslant i\leqslant p^{q_i}$, does not pass through the edges connecting neighboring copies of the graph $\Lambda_{i, r_i-1}$, hence its length is the same as in the graph $\Lambda_{i, r_i-1}$.

Let's prove now that for each $i, m+1\leqslant i\leqslant n, i\neq h\  |\varphi_{\Lambda}(u_i)|\geqslant|\varphi_{\Lambda}(u_1)|$. When we construct graphs $\Lambda_{1i}$ thanks to the condition on $h$ it is true that for each $i, 1\leqslant i\leqslant r_1-1,$ and for each $s, n+1>s>m$ the path $S_1^i$ is contained in all maximal $u_s$-cycles. That is when $1\leqslant i<r_1-1$ after we construct $\Lambda_{1, i+1}$ from $\Lambda_{1i}$ the length of the maximal $u_s$-cycle increases in $p$ times. When we step to $\Lambda_{1, r_1}$ from $\Lambda_{1, r_1-1}$ the length of the maximal $u_s$-cycle increases either in $p^{q_1-1}$ times or in $p^{q_1}$ times. So for all $s, n+1>s>m, s\neq h\
|\varphi_{\Lambda}(u_s)|\geqslant|\varphi_{\Lambda}(u_1)|$. Thus it is proved that for pairwise incommensurable elements $u_1, ..., u_n$ and for homomorphism $\varphi_m$ which satisfies the condition that $|\varphi_m(u_1)|=...=|\varphi_m(u_m)|<|\varphi_m(u_i)|, i>m$, and $u_j$- and $v$-cycles are simple in $Cay(\varphi_m(G); \varphi_m(Z)), 1\leqslant j\leqslant n$ there exists a homomorphism $\varphi_{m+1}'$ such that $|\varphi_{m+1}'(u_1)|=...=|\varphi_{m+1}'(u_m)|\leqslant|\varphi_{m+1}'(u_i)|, n+1>i>m$ and for some $h, m<h<n,$ which does not depend on $\varphi_m', |\varphi_m'(u_h)|/|\varphi_m'(u_1)|<|\varphi_m(u_h)|/|\varphi_m(u_1)|$. Besides $u_j$- and $v$-cycles are simple in $Cay(\varphi_{m+1}'(G); \varphi_{m+1}'(Z)), 1\leqslant j\leqslant n$. Repeating this procedure and changing indices we can construct the homomorphism $\varphi_{m+1}$ such that $|\varphi_{m+1}(u_1)|=...|\varphi_{m+1}(u_m)|=|\varphi_{m+1}(u_{m+1})|\leqslant|\varphi_{m+1}(u_i)|,
n+1>i>m+1$. Besides $u_i$- and $v$-cycles are simple in $Cay(\varphi_{m+1}; \varphi_{m+1}(Z))$. Repeating this process we shall come to the homomorphism $\varphi_n$ of $F$ onto a finite $p$-group such that the images of elements $u_1,..., u_n$ pairwise coincide and all $u_i$- and $v$-cyles are simple in the graph $Cay(\varphi_n(F); \varphi_n(Z)), 1\leqslant i\leqslant n$. Suppose that $|\varphi_n(u_1)|\leqslant|\varphi_n(v)|$. In this case same as at the construction of the homomorphism $\varphi_{m+1}'$ from the homomorphism $\varphi_m$ we are able to construct the homomorphism $\varphi_n'$ of $F$ onto a finite $p$-group such that $|\varphi_n'(u_1)|=...=|\varphi_n'(u_n)|$, and $|\varphi_n'(u_1)|/|\varphi_n'(v)|>|\varphi_n(u_1)|/|\varphi_n(v)|$ and $u_j$- and $v$-cycles are simple in $Cay(\varphi_n'(F); \varphi_i'(Z))$.
Applying this procedure the appropriate number of times we shall obtain the required homomorphism $\psi$. It easy to notice that when we step from the homomorphism $\varphi_m$ to the homomorphism $\varphi_m'$ the orders of images of elements $u_i, 1\leqslant i\leqslant k$ can only increase. Hence $|\psi(u_1)|>N$. Lemma 1 is proved.

\textbf{Lemma 2.} \textsl{Let $G=F_1{\ast\atop a=b}F_2$, where $F_1, F_2$ are free groups, $a\in F_1\setminus\{1\}, b\in F_2\setminus\{1\}$ are not proper powers. Consider the elements $u=a_1b_1...a_nb_n, v=a_1'b_1'...a_m'b_m'\in G, a_i, a_j'\in F_1\setminus\<a\>, b_i, b_j'\in F_2\setminus\<b\>, 1\leqslant i\leqslant n, 1\leqslant j\leqslant m$. For each element $q\in\<a\>\setminus\{1\}$ there exist homomorphisms $\varphi, \psi$ of groups $F_1$ and $F_2$ onto finite groups such that
$|\varphi(a)|=|\psi(b)|>1$. Besides the nonidentity elements from the set
$\Delta_{F_1}$:
$$
\Delta_{F_1}=\{a_1^{-1}qa_1, a_i, a_j', a_ia_k^{-1}, a_j'a_l'^{-1},
a_ia_j'^{-1}| 1\leqslant i, k\leqslant n, 1\leqslant j, l\leqslant
m\}
$$
have nonidentity images after $\varphi$. Elements of the set $\Delta_{F_1}$ which do not belong to $\<a\>$ possess the images which also do not belong to $\<\varphi(a)\>$. The analogical conditions is true for the homomorphism $\psi$ and subgroup $\<b\>$ and for the set $\Delta_{F_2}$:
$$
\Delta_{F_2}=\{b_i, b_j', b_ib_k^{-1}, b_j'b_l'^{-1}, b_ib_j'^{-1}|
1\leqslant i, k\leqslant n, 1\leqslant j, l\leqslant m\}
$$
}

\textbf{Proof.}

Since subgroups $\<a\>, \<b\>$ are $p'$-isolated for each prime number $p$ then the groups $F_1, F_2$ residually $p$-finite with respect to entry in subgroups $\<a\>, \<b\>$ and there exist homomorphisms $\varphi_1, \psi_1$ of $F_1$ and $F_2$ correspondingly onto finite $p$-groups satisfying the following properties. The elements of the set $\Delta_{F_1}$ which do not belong to the subgroup $\<a\>$ have images which do not belong to the subgroup $\varphi_1(\<a\>)$ and nonidentity elements of the set $\Delta_{F_1}$ have nonidentity images. The analogical properties are true for the homomorphism $\psi_1$ for the set $\Delta_{F_2}$ and for the subgroup $\<b\>$. Let prime number $p$ is greater than the maximum of lengths of elements $a$ and $b$, and $N=$ max $(|\varphi_1(a)|, |\psi_1(b)|)$. According to lemma 1 there exists the homomorphism $\epsilon$ of $F_1\ast F_2$ onto a finite $p$-group which is so that $|\varphi_2(a)|=|\psi_2(b)|>N$. Put $\varphi=\epsilon_{F_1}, \psi=\epsilon_{F_2}$. It is easy to notice that the homomorphisms $\varphi: F_1\rightarrow \varphi_1(F_1)\times\varphi_2(F_1), \varphi: F_2\rightarrow \psi_1(F_2)\times\psi_2(F_2)$ defined by formulas $\varphi: g\mapsto \varphi_1(g)\varphi_2(g), \psi: f\mapsto \psi_1(f)\psi_2(f), g\in F_1, f\in F_2$ are as required.
Lemma 2 is proved.

Now we need to define the following notation. Consider edges $e$ and $f$ and paths $E=\{e_i| i\in I\}, F=\{f_i| i\in I\}$ of an action graph of the group $A{\ast\atop C}B$. We shall say that $e$ and $f$ are $C$-near if $C(\alpha(e))=C(\alpha(f))$ and
$C(\omega(e))=C(\omega(f))$ and labels of edges $e, f$ belong to one free factor. We shall also say that the paths $E$ and $F$ are $C$-near if for each $i\in I$ the edges $e_i$ and $f_i$ are $C$-near.

\section{Proof of theorem 1.}

Put $G=A{\ast\atop C=D}B$ where $A, B$ are free, $C=\<a\>\subset A, D=\<b\>\subset B$. Consider elements $u, v\in G$ so that $u$ does not conjugate to $v^{\pm1}$. We shall consider that $u$ and $v$ are cyclically reduced.

Consider the case when $u, v\in A$. In [1] it was proved that there exists the homomorphism $\varphi_1$ of $A$ onto a finite group such that $\mid\varphi_1(u)\mid\neq\mid\varphi_1(v)\mid$. Put $\mid\varphi_1(a)\mid=n$. According to theorem 4 there exists the homomorphism $\varphi_2$ of $F_2$ onto a finite group such that $\mid\varphi_2(b)\mid=n$. Since the map $\rho:
\varphi_1(a)\mapsto \varphi_2(b)$ is completed till isomorphism between the subgroups $\<\varphi_1(a)\>_n$ and $\<\varphi_2(b)\>_n$, there exists exactly one homomorphism $\varphi_1\ast\varphi_2$ of the group $G$ onto the group $G_1=\varphi_1(A){\ast\atop\varphi_1(C)\simeq\varphi_2(D)}\varphi_2(B)$ such that $\varphi_1\ast\varphi_2\mid_A=\varphi_1, \varphi_1\ast\varphi_2\mid_B=\varphi_2$. Since the group $G_1$ is residually finite and the elements $u_1=\varphi_1(u)$ and
$v_1=\varphi_2(v)$ have different orders there exists a homomorphism $\varphi_3$ of $G_1$ onto a finite group such that $\<\varphi_1(u)\>\simeq\<\varphi_3(u_1)\>, \<\varphi_1(v)\>\simeq\<\varphi_3(v_1)\>$.

Suppose now that $u\in A\setminus\{C\}, v\in B\setminus\{C\}$. Consider the free group $F=A\ast B$. According to lemma 1 there exists the homomorphism $\varphi_1$ of $F$ onto a finite group such that $\mid\varphi_1(u)\mid=\mid\varphi_1(a)\mid=\mid\varphi_1(b)\mid>\mid\varphi_1(v)\mid$. Let $\psi_1=\varphi_1|_A, \psi_2=\varphi_1|_B$. Since $|\psi_1(a)|=|\psi_2(b)|$ there exists exactly one homomorphism $\psi_1\ast\psi_2$ of $A\ast_C B$ onto $\psi_1(A)\ast_{\psi_1(C)}\psi_2(B)$ such that $\psi_1\ast\psi_2|_A=\psi_1, \psi_1\ast\psi_2|_B=\psi_2$. Besides $|\psi_1\ast\psi_2(u)|>|\psi_1\ast\psi_2(v)|$. The group $\psi_1(A)\ast_{\psi_1(C)}\psi_2(B)$ is residually finite so there exists a homomorphism of this group onto a finite group such that the images of elements $\psi_1\ast\psi_2(u), \psi_1\ast\psi_2(v)$ have different orders.

In order to consider the most common case let's change notations slightly.

Consider the group $\widetilde{G}=\widetilde{A}{\ast\atop\widetilde{C}=\widetilde{D}}\widetilde{B}, \widetilde{A}, \widetilde{B}$ are free, $\widetilde{C}=\<\widetilde{a}\>, \widetilde{D}=\<\widetilde{b}\>, \widetilde{a}$ and $\widetilde{b}$ are not proper powers. Consider elements $\widetilde{u}, \widetilde{v}$, $\widetilde{u}\notin \widetilde{A}\cup \widetilde{B}$, $\widetilde{u}$ and $\widetilde{v}$ are cyclically reduced. Also $\widetilde{u}$ is not conjugate to $\widetilde{v}^{\pm1}$. Fix some reduced notations for the elements $\widetilde{u}$ and $\widetilde{v}$. We may consider that they have the following presentations: $\widetilde{u}=\widetilde{a_1}\widetilde{b_1}...\widetilde{a_m}\widetilde{b_m}, \widetilde{v}=\widetilde{a_1}'\widetilde{b_1}'...\widetilde{a_l}'\widetilde{b_l}'$ or $\widetilde{v}\in A\cup B$ and in this case we may consider that $\widetilde{v}=\widetilde{a_1}', l=1, \widetilde{a_i}, \widetilde{a_j}'\in \widetilde{A}\setminus\widetilde{C}, \widetilde{b_i}, \widetilde{b_j}'\in \widetilde{B}\setminus\widetilde{D}, 1\leqslant i\leqslant m, 1\leqslant j\leqslant l$. According to lemma 2 there exist homomorphisms $\varepsilon, \sigma$ of groups $\widetilde{A}, \widetilde{B}$ onto finite groups $A$ and $B$ such that $|\varepsilon(\widetilde{a})|=|\sigma(\widetilde{b})|$ and $|\varepsilon(\Delta_A\setminus\{1\})|=|\Delta_A\setminus\{1\}|, |\varepsilon(\Delta_A\setminus\{\widetilde{A}\})|=|\Delta_A\setminus\{\widetilde{A}\}|, |\sigma(\Delta_B\setminus\{1\})|=|\Delta_B\setminus\{1\}|, |\sigma(\Delta_B\setminus\{\widetilde{B}\})|=|\Delta_B\setminus\{\widetilde{B}\}|$. If $\widetilde{v}\notin A\cup B$ then $\Delta_A=\{\widetilde{a_1}^{-1}\widetilde{q}\widetilde{a_1}, \widetilde{a_i}, \widetilde{a_j}', \widetilde{a_i}\widetilde{a_k}^{-1}, \widetilde{a_j}'\widetilde{a_t}'^{-1}, \widetilde{a_i}\widetilde{a_j}'^{-1}| 1\leqslant i, k\leqslant m, 1\leqslant j, t\leqslant l\}, \Delta_B=\{\widetilde{b_i}, \widetilde{b_j}', \widetilde{b_i}\widetilde{b_k}^{-1}, \widetilde{b_j}'\widetilde{b_t}'^{-1}, \widetilde{b_i}\widetilde{b_j}'^{-1}| 1\leqslant i, k\leqslant m, 1\leqslant j, t\leqslant l\}$ and $\widetilde{q}$ equals unit if elements $\widetilde{u}$ and $\widetilde{v}$ do not belong to one left coset on the amalgamated subgroup, otherwise it is the element of the amalgamated subgroup of the group $\widetilde{G}$ such that $\widetilde{v}=\widetilde{u}\widetilde{q}$. If $\widetilde{v}\in\widetilde{A}$ then $\Delta_A=\{\widetilde{a_i}, \widetilde{a_1}'| 1\leqslant i\leqslant m\}, \Delta_B=\{\widetilde{b_i}| 1\leqslant i\leqslant m\}$. Thus there exists exactly one homomorphism $\rho$ of $\widetilde{G}$ onto a group $\varepsilon(\widetilde{A}){\ast\atop\varepsilon(\widetilde{C})}\sigma(\widetilde{B})$, such that $\rho|_{\widetilde{A}}=\varepsilon, \rho|_{\widetilde{B}}=\sigma$ besides according to the conjugacy theorem for free products with amalgamation [7] $\rho(\widetilde{u})$ does not conjugate to $\rho(\widetilde{v})$, $\rho(\widetilde{u})$ does not conjugate to $\rho(\widetilde{v}^{-1})$ and $\rho(\tilde{u})$ does not belong to a subgroup which is conjugate to either $\rho(\tilde{A})$ or $\rho(\tilde{B})$. Thereby we may suppose that free factors are finite. Put $\rho(\widetilde{G})=G, \rho(\widetilde{u})=u, \rho(\widetilde{v})=v, \rho(\widetilde{C})=C, \rho(\widetilde{a_i})=a_i, \rho(\widetilde{b_i})=b_i, 1\leqslant i\leqslant m$. According to theorem 3 there exists a homomorphism $\varphi$ of $G$ onto a finite group such that in the Cayley graph of the group $\varphi(G)\ \Gamma=Cay(\varphi(G); \{\varphi(A\cup B)\})$ the representatives of $u$-cycles have no $1$-near vertices and if $v\notin A\cup B$ the representatives of $v$-cycles in this graph also have no $1$-near vertices. Suppose that $|\varphi(u)|=|\varphi(v)|$. Consider the inductive process of construction of action graphs of the group $G$ $\Gamma_k, 1\leqslant k\leqslant s$, which satisfy the following conditions ($s$ is either the nonnegative integer or the infinite symbol):

1) The length of each $u$-cycle divides the length of the maximal $u$-cycle. The same is true for $v$-cycles.

2) All representatives of $u$-cycles have no $1$-near vertices. The analogical condition is true for $v$-cycles in case $v\notin A\cup B$.

3) When $k\leqslant s$ there exist two paths $S_1^k$ and $S_2^k$ of length $k$, so that $S_1^k$ is contained in some representative of maximal $u$-cycle and all representatives of all maximal $v$-cycles have subpaths which are $C$-near to the path $S_2^k$. Besides $S_1^k$ and $S_2^k$ are $C$-near.

4) The lengths all maximal $u$- and $v$-cycles coincide.

Representatives of maximal $u$ and $v$-cycles will be referred to as maximal representatives of $u$- and $v$-cycles.

If $v\in A\cup B$ we use the residual finiteness of the group $G$ and that the order of $u$ is infinite while the order of $v$ is finite. Let's construct the graph $\Gamma_1$. Put $n=|\varphi(u)|$. Fix the $u$-cycle in the graph $\Gamma$ going from some vertex $p$ and fix also its representative $S=e_1...e_m$ whose label equals $u^n$. Let $q=\omega(e_1)$. Consider $n$ copies of the graph $\Gamma: \Delta_1,..., \Delta_n$. From each copy we delete edges with labels from $A\setminus C$ whose begin or end points belong to the subgraph $C(q_i)$. For each $i, 1\leqslant i\leqslant n$ (indices are modulo $n$) and for each edge that was deleted we shall make the following transformation. Consider the edge $f$ of the graph $\Gamma$ whose terminal point $q'$ belongs to $C(q)$ and whose label from $A\setminus C$. Another terminal point of this edge is denoted as $p'$. In graphs $\Delta_i$ and $\Delta_{i+1}$ connect vertices $p'^i$ and $q'^{i+1}$ by the new edge $f_i$ whose label coincide with the label of the edge $f$. If the edge $f$ goes away from the vertex $p'$ then the edge $f_i$ goes away from the vertex $p'^i$. If $f$ goes into $p'$ then $f_i$ goes into $p'^i$. The obtained graph is denoted as $\Gamma_1'$. This graph satisfies the condition 1) because of the condition on the number of copies of the graph $\Gamma$ and on the value of $|\varphi(u)|$. The condition 2) is also held for $\Gamma_1'$ because otherwise it is not true for the graph $\Gamma$. Consider the representative of the $u$-cycle in $\Gamma_1'$ which goes away from the vertex $w=\alpha(e_1)^1$ and which corresponds to the same notation of $u$ as $S$. Its second edge is the edge $f=e_2^2$. Consider $n^2$ copies of the graph $\Gamma_1': \Lambda_1,...,
\Lambda_{n^2}$. From each graph $\Lambda_i$ we delete all edges whose terminal points belong to the subgraph $C(\omega(f^i))$, and whose labels belong to $B\setminus C$. Consider each edge $h$ of the graph $\Gamma_1'$ whose terminal point $q''$ belong to $C(\omega(f_i))$ and whose label is in $B\setminus C$. The second terminal point of $h$ is the vertex $p''$. For each $i, 1\leqslant i\leqslant n^2$ (indices are modulo $n^2$) connect vertices $p''^i$ and $q''^{i+1}$ of graphs $\Lambda_i$ and $\Lambda_{i+1}$ by the edge $g_i$. If $h$ goes into $p''$ then $g_i$ goes into $p''^i$ and if $h$ goes away from $p''$ then $g_i$ goes away from $p''^i$, Lab $(h)=$ Lab $(g_i)$. The obtained graph will be denoted as $\Gamma_1$ and $h_1$ is the new edge in $\Gamma_1$ which connect vertices $\alpha(f)^1$ and $\omega(f)^2$ and which was obtained after we deleted the edge $f^1$. The correctness of properties 1), 2) for $\Gamma_1$ is proved analogically as for the graph $\Gamma_1'$. Furthermore by the construction all representatives of all maximal $v$-cycles of the graph $\Gamma_1'$ pass through the vertex from $C(\alpha(f))$. Thus according also to the condition about the absence of $1$-near vertices for the representatives of $u$- and $v$-cycles we deduce that one of the following two alternatives is held. Let $p$ be the vertex of the graph $\Gamma_1'$ which is the beginning of maximal $v$-cycle. That is the representatives of this $u$-cycle pass through vertices from $C(\alpha(f))$. Then either the corresponding vertex $p^1$ of the graph $\Gamma_1$ is the beginning of the $v$-cycle whose representatives pass through edges which are $C$-near to the edge $h_1$ or the representatives of each such $v$-cycle pass only through vertices of the singular copy of the graph $\Gamma_1'$ because the representatives of these $v$-cycles pass through edges whose labels belong to $B\setminus C$ and their terminal vertices belong to $C(\alpha(h_1))$ but these edges are not $C$-near for the edge $h_1$ and since representatives of $v$-cycles in $\Gamma_1$ do not have $1$-near vertices the representatives of such $v$-cycle do not pass through edges which are $C$-near to $h_1$ and does not contain vertices from $C(\omega(h_1))$. In the last case each $v$-cycle in the graph $\Gamma_1$ has length which is less than the length of the maximal $u$-cycle in the same graph. Otherwise we have the graph $\Gamma_1$ which satisfies the properties 1), 2), 3), 4) if we put $S_1^1=h_1$ and $S_2^1$ is an arbitrary edge which is $C$-near to the edge $h_1$.

Suppose we have the graph $\Gamma_k, k\geqslant1$. Let $n$ be the length of the maximal $u$-cycle in $\Gamma_k$. Let $T$ be the representative of the maximal $u$-cycle passing through $S_1^k$, $T=e_1...e_l$. Then $S_1^k=e_{q+1}...e_{q+k}$. Consider $n$ copies of the graph $\Gamma_k: \Delta_1,..., \Delta_n$. Let $q$ be the end of the edge $e_{q+k+1}$ in the graph $\Gamma_k$. From each graph $\Delta_i$ delete all edges whose labels belong to the same free factor as the label of the edge $e_{q+k+1}^i$ but do not belong to $C$ and whose terminal vertices belong to $C(q_i)$. Further for each $i, 1\leqslant
i\leqslant n$ (indices are modulo $n$) and for each edge that was deleted the following procedure will be made. Let $f$ be a positively oriented edge of $\Gamma_k$ whose copy was deleted in $\Delta_i$,
$r=\alpha(f), s=\omega(f)$. Let $s'$ be the vertex from $\{r; s\}\cap C(q)$, and $r'$ is the vertex from $\{r; s\}\setminus\{s'\}$. Connect vertices $r'^i$ and $s'^{i+1}$ of graphs $\Delta_i, \Delta_{i+1}$ correspondingly by the edge $h_i$. Its label equals the label of the edge $f$. If $f$ goes into the vertex $r'$ then $h_i$ goes into $r'^i$. If $f$ goes away from the vertex $r'$ then $h_i$ goes away from $r'^i$. The new graph will be referred to as $\Gamma_{k+1}'$. The new graph satisfies the properties 1), 2). The proof of this fact is analogical to the case of the graph $\Gamma_1$. Put $S_2^k=\varepsilon_1...\varepsilon_k$ where the edge $\varepsilon_i$ is $C$-near to the edge $e_{q+i}, 1\leqslant i\leqslant k$. Let $h_1$ be the new edge in $\Gamma_{k+1}'$ that connect vertices $\alpha(e_{q+k+1})^1$ and $\omega(e_{q+k+1})^2$ and which was appended instead of the edge $e_{q+k+1}^1$. Fix the representative of some $v$-cycle in $\Gamma_k$ which has the subpath $C$-near to the path $S_2^k$. Suppose this subpath goes away from the vertex $w$. Consider the corresponding representative $Q$ in $\Gamma_{k+1}'$ going away from the vertex $w_1$. It has the subpath which is $C$-near to the path $(S_2^k)^1$. Consider the edge $d$ of $Q$ which is the next after the edge $C$-near to the edge $\varepsilon_k^1$. If $d$ is $C$-near to the edge $h_1$ then the length of the $v$-cycle $W$ corresponding to $Q$ coincides with the length of the $u$-cycle maximal $u$-cycle of the graph $\Gamma_{k+1}'$ and $v$-cycle $W$ is the maximal in the graph. In this case the property 4) is held for $\Gamma_{k+1}'$. The property 3) fulfilled since we may put $S_1^{k+1}=(S_1^k)^1\cup h_1$ and $S_2^{k+1}=(S_2^k)^1\cup z$ where $z$ is the edge which is $C$-near to the edge $d$ and incident to the edge $\varepsilon_k^1$. Thus $\Gamma_{k+1}'=\Gamma_{k+1}$. Note that there is no edge $d$ whose terminal point belongs to $C(\alpha(h_1))\cup C(\omega(h_1))$ and which belongs to $Q$ and which satisfies the following condition. The edge $d$ is the next edge is not the next edge after the edge of $Q$ which is $C$-near to the edge $\varepsilon_k^1$ because in this case we have 1-near vertices in $Q$. If for some maximal representative of the $v$-cycle the edge $d$ defined above is not $C$-near to the edge $h_1$ then the length of the corresponding $v$-cycle is less than the length of the maximal $u$-cycle in $\Gamma_{k+1}'$ and it has empty intersection with $\Delta_i, i>1$. Furthermore for each representative of the $v$-cycle of $\Gamma_{k+1}'$ which goes away from some vertex $w^j$ and pass through the subgraphs $\Delta_i$ for different $i$ the corresponding representative of the  $v$-cycle in $\Gamma_k$ going away from the vertex $w$ is not maximal in the graph $\Gamma_k$. So in this case the length of the maximal $u$-cycle is greater than the length of the maximal $v$-cycle in $\Gamma_{k+1}'$ and theorem is proved. So we consider that for each natural $k$ there exists the graph $\Gamma_k$. It means that in the original Cayley graph $\Gamma$ the representatives of $u$- and $v$-cycles going from the same vertex are $C$-near. Since the orders of images of $u$ and $v$ are equal and representatives of $u$- and $v$-cycles have no 1-near vertices in $\Gamma$ the lengths $u$ and $v$ coincide. Fix some irreducible notations for $u$ and $v$: $u=a_1b_1...a_mb_m, v=a_1'b_1'...a_m'b_m'$. Consider representatives $Q$ and $R$ of $u$- and $v$-cycles going from one vertex and corresponding to the chosen notations. Since $Q$ and $R$ are $C$-near then $a_1$ and $a_1'$ belong to one coset on $C$ so we may consider that $a_1=a_1'$. Using this argument we may also consider that $a_i=a_i', b_j=b_j', 1\leqslant i\leqslant m,
1\leqslant j\leqslant m-1, b_m=b_m'h$ for some $h\in C$. That is $u$ and $v$ belong to one coset on $C$: $u=a_1b_1...a_mb_m, v=a_1b_1...a_mb_mh$. Since $Q$ and $R$ are $C$-near then $\varphi(a_1^{-1}ha_1)\in \varphi(C)$. Due to the definition of the irreducible notations [8] we have $\widetilde{a_1}\notin\widetilde{C}$. Since $\widetilde{a}$ is not a proper power we conclude that $\widetilde{a_1}^{-1}\widetilde{q}\widetilde{a_1}\notin\widetilde{C}$. Furthermore $\widetilde{v}=\widetilde{u}\widetilde{q}$ for some  $\widetilde{q}\in\<\widetilde{C}\>$. But this contradicts the conditions on $\rho$.

The case when $v=1$ is trivial. Theorem 1 is proved.

\begin{center}
\large{Acknowledgements.}
\end{center}

The author thanks A. A. Klyachko for valuable comments and discussions.

\begin{center}
\large{References.}
\end{center}

1. \textsl{Klyachko A. A.} Equations over groups, quasivarieties,
and a residual property of a free group // J. Group Theory. 1999.
\textbf{2}. 319--327.

2. \textsl{Wise, Daniel T.} Subgroup separability of graphs of free groups with cyclic edge groups.
Q. J. Math. 51, No.1, 107-129 (2000). [ISSN 0033-5606; ISSN 1464-3847]

3. \textsl{Graham A. Niblo.} H.N.N. extensions of a free group by \textbf{z} which are subgroup separable. Proc. London Math. Soc. (3), 61(1):18--23, 1990.

4. \textsl{Yedynak V. V.} Separability with respect to order. \textsl{Vestnik
Mosk. Univ. Ser. I Mat. Mekh.} \textbf{3}: 56-58, 2006.

5. \textsl{Romanovskii N. S.} On the residual finiteness of free products with respect to entry // Izv. AN SSSR. Ser.
matem., 1969, \textbf{33}, 1324-1329.

6. \textsl{P. Stebe.} Conjugacy separability of certain free
products with amalgamation. Trans. Amer. Math. Soc. 156 (1971).
119--129.

7. \textsl{Kargapolov, M. I., Merzlyakov, Yu. I.} (1977). Foundations of group
theory. \textsl{Nauka}.

8. \textsl{Lyndon, R. C., Schupp, P. E.} (1977). Combinatorial group
theory.\ \textsl{Springer-Verlag}.

\end{document}